# New result in NG-Groups

Faraj.A.Abdunabi[1], Ahmed shletiet[2]

[1]Lecturer in Mathematics Department, AJdabyia University

Faraj.A.Abdunabi@uoa.edu.ly

[2]Lecturer in Mathematics Department, AJdabyia University

shletiet@ uoa.edu.ly

**Abstract**: We investigate and present a new result in NG groups that consisting of consisting of non-bijective transformations cannot be subset of symmetric groups. In this paper, we revive the concepts of *NG* groups. Moreover, we introduce the definition of regular NG and inverse NG the new concept as new results in these groups. Hover, some properties are studied.

1- **Introduction:**

The transformation group is important and good part of the group theory [1]. We recall a permutation group on a non-empty set *A* is a group consisting of bijections transformations from *A* to itself with respect to compositions mapping on anon empty set A which the one most important transformation groups for more details see [2],[3],[4]. It is well known that any permutation group on a set *A* with cardinality *n* has order not greater that *n*!. In [5], Y. Wu, X.Wei present the conditions of the groups generated by nonobjective transformations on a set. Author in [6], introduce find the maximal order of these groups. In this paper, we fined a new result in these groups, which is the inverse of these groups.

2- **Preliminary :**

We review some basics definitions and properties of the finite group theory that we will used in our paper. For more detailed in lots of abstract algebra and finite group theory can see [7],[8] would be good supplementary sources for the theory needed here.

**Definition 2.1**. Suppose that *A* be non-empty set. A binary relation *R* in *A* is called *an equivalence relation* on *A* if it satisfy the following*:*
   *(i)*     *xRx* for any *x*∈*A*;
   *(ii)*    for any *x,y*∈ *A*, if *xRy* then *yRx* ;
   *(iii)*   for any *x,y,z*∈*A*, if *xRy* and *yRz* then *xRz.*

We denoted *P(A)* to the set of all its transforms. For any *f*∈*P(A)*, we use Im(*f*) to denoted the image of *f*. Also, *Z* and $Z_>$ will respective dented the set of integers and positive integers.

**Definition 2.2**. Suppose that $R$ is e an equivalence relation on $A$. For an element $x \in A$, we call $\{x \in A / x R\}$ the equivalence class of a determined by $R$, which is denoted by $R[x]$. And $A/R = \{R[x] | x \in A\}$ is called the quotient set of $A$ relative to the equivalence relation $R$.

**Lemma 2.3**. *([5],theorem 1)* *For any $f \in NG$ and the e the identity element of NG, $R_e = R_f$.*
*Proof.*

For any $x \in A$, our goal is to show that $[a]_f = [a]_e$. On one hand, if $a \in [x]_f$, i.e. $f(a) = f(x)$. Since $NG$ is a group with identity element $e$, there is a transformation $f' \in NG$ such that $f'f = e = ff'$. Therefore, $e(a) = f'(f(a)) = f'(f(x)) = e(x)$, Which yields that $a \in [x]_e$. On the other hand, if $y \in [a]_e$ i.e. $e(a)e(y)$. Hence, $f(a) = (fe)(a) = f(e(y)) = (fe)(y) = f(y)$, which implies $y \in [x]_f$. It follows that $[x]_e = [x]_f$ for any $x \in A$, as wanted.

**Remark 2.4.** For Lemma 2.3, we see that $R_f = R_g$ for any element $f, g \in NG$. The following Theorem is revised version of Theorem 2,[5].

**Theorem 2.5.** *Let $f$ be an element in $P(A)$ and $\hat{f}$ be the induced transformation of $f$ on $A/R_f$, i.e $\hat{f}: A/R_f \longrightarrow A/R_f$, $[x]_f \mapsto [f(x)]_f$. Then the following hold:*
*(i)    There exists a groups $NG \subseteq P(A)$ containing $f$ as the identity element iff $f^2 = f$.*
*(ii)    There is a groups $NG \subseteq P(A)$ containing $f$ as the identity element iff $\hat{f}$ is a bijective on $A/R_f$.*

We make some corrections to the original proofs of corollaries are from [5]. Actually, we adopt the restriction of finiteness on $A$ in the first corollary from the original one. And we used the finiteness on $A$ in the second corollary; the original one did not use it.

**Corollary 2.6**. *Let $f$ be an element in $P(A)$. Then $f^2 = f$ iff the induced mapping $\hat{f}$ on $A/R_f$ is the identity element.*
*Proof.*

($\Rightarrow$) we suppose that $f^2 = f$. Then for any $[x]_f \in A/R_f$, as $f(x) = f(f(x))$ we see that $[x]_f = [f(x)]_f$. It follows that $\hat{f}([x]_f) = [f(x)]_f = [x]_f$; which implies that $\hat{f}$ is the identity mapping on $A/R_f$.

($\Leftarrow$) we assume that $\hat{f}$ is the identity mapping on $A/R_f$. Then for any $[x]_f \in A/R_f$, the condition that $\hat{f}([x]_f) = [x]_f$ will imply that $[f(x)]_f = [x]_f$ and hence $f(f(x)) = f(x)$. It follows that $f^2 = f$ as required.

**Corollary 2.7.** *Suppose that $A$ is a finite set and $f$ is an element in $P(A)$. Then there is a group $NG \subseteq P(A)$ containing $f$ as an element iff $Im(f) = Im(f^2)$.*

*Proof.*
($\Rightarrow$) we suppose that there is a group $NG \subseteq P(A)$ containing $f$ as an element. Let $e$ be the identity element of $NG$. Then by Theorem 2.5, the induced mapping $\hat{f}$ is a bijection on $A/R_f$. In particular, $\hat{f}$ is surjective and thus for any $x \in A$, there is a $[y]_f \in A/R_f$ such that $\hat{f}([y]_f) = [x]_f = [f(y)]_f$; which yields that $f(x) = f(f(y)) = (f^2)(y)$. As a result, $Im(f) \subseteq Im(f^2)$ and thus $Im(f) = Im(f^2)$.

($\Leftarrow$) we suppose that Im($f$) = Im($f^2$). Thus, for any $f(x) \in$ Im($f$) there is a $y \in A$ such that $f(x) = f(f(y))$ and hence $\hat{f}([y]f) = [x]f$; which implies that $\hat{f}$ is surjective on $A/R_f$. Note that we are assuming that $A$ is finite and so is $A/R_f$. We have that the induced mapping $\hat{f}$ is bijective. By Theorem 2.5, the assertion follows.

**Remark 2.2.** Suppose that $NG \subseteq P(A)$ is a group. We have seen, in Remark 2.1, that $R_f = R_g$ for any elements in $NG$ and we will denote the common equivalence relation by $R$. Also, by Theorem 2.5, each element $f \in NG$ will induce a bijection $\hat{f}$ on $A/R$.

The following theorem is crucial since it turns a group $NG \subseteq P(A)$ into a permutation group.

**Theorem 2.8.** Suppose that $NG \subseteq P(A)$ is a group. Set $\hat{NG} = \{\hat{f} | f \in NG\}$; then $\hat{NG}$ is a permutation group on A/R and $\rho : NG \to \hat{NG}, f \mapsto \hat{f}$, is an isomorphism.

*Proof.*

For any $f, g \in NG$ and any $[x] \in A/R$, we have $\rho(fg)([x]) = [(fg)(x)] = [f(g(x))] = \rho(f)([g(x)]) = (\rho(f)\rho(g))([x])$; which implies that $\rho(fg) = \rho(f)\rho(g)$ and thus $\rho$ is a homomorphism. By the definition of $\hat{NG}$, it is obvious that $\rho$ is surjective.

Now suppose that $\rho(f) = \rho(g)$ for two elements $f, g \in NG$, i.e. $[f(x)] = [g(x)], \forall x \in A$: Let $e$ be the identity element of $NG$, then we have $[f(x)]_e = [g(x)]_e; \forall x \in A$. It follows that $e(f(x)) = e(g(x)); \forall x \in A$. Hence, $f(x) = (ef)(x) = e(f(x)) = e(g(x)) = g(x), \forall x \in A$, and therefore $f = g$. We conclude that $\rho$ is injective. As a consequence, $\rho$ is an isomorphism.

**Examples 2-9** Let $F$ be a field, and $V$ be a vector space of dimension 2 over $F$. And $\{v_1, v_2\}$ a basis of $V$. For any $a \in F^*$, we define $T_a : V \to V$, $\mu = a_1 v_1 + a_2 v_2 \to a\, a_1 v_1$.

In other words, the linear transformation $T_a$ has matrix $\begin{pmatrix} 1 & 0 \\ 0 & 0 \end{pmatrix}$ with respect to the ordered basis $v_1, v_2$

Set $NG = \{T_a : a \in F^*\}$, then $NG$ is group. But, $NG$ is not a subset of symmetric group.

**Examples 2-10** Suppose $X = \{1,2,3\}$. We know, $S_3 = \{(1,2,3), (2,3,1), (3,1,2), (1,3,2), (3,2,1), (2,1,3)\}$ is a group but it is not an abelian group. We introduce the elements of Trans(X) as: $\{(1,1,1), (1,1,2), (1,1,3), (1,2,1), (1,2,2), (1,2,3), (1,3,1), (1,3,2), (1,3,3), (2,1,1), (2,1,2), (2,1,3), (2,2,1), (2,2,2), (2,2,3), (2,3,1), (2,3,2), (2,3,3), (3,1,1), (3,1,2), (3,1,3), (3,2,1), (3,2,2), (3,2,3), (3,3,1), (3,3,2), (3,3,3)\}$. There exits some groups that are subsets of *Trans(X)*, but not subsets of the $S_n$.

We make composition of all transformations we get:

If we apply same way, we get

The groups of order 2 are : $G_1 = \{(1,1,3),(3,3,1)\}$, $G_2 = \{(1,2,1),(2,1,2)\}$,
$G_3 = \{(1,2,2),(2,1,1)\}$,
$G_4 = \{(1,2,3),(1,3,2)\}$, $G_5 = \{(1,2,3),(2,1,3)\}$, $G_6 = \{(1,2,3),(3,2,1)\}$, $G_7 = \{(1,3,3),(3,1,1)\}$,
$G_8 = \{(2,2,3),(3,3,2)\}$, and $G_9 = \{(2,3,2),(3,2,3)\}$. But, $G_4 = \{(1,2,3),(1,3,2)\}$,
$G_5 = \{(1,2,3),(2,1,3)\}$, and $G=6 = \{(1,2,3),(3,2,1)\}$ are subsets of Sym(X).

**Proposition 2-11:** Let $NG_1$ and $NG_2$ be two NG-groups, then $NG_1 \cup NG_2$ and $NG_1 \cap NG_2$ are not necessary to be *NG*-groups.

*Proof:*
Let X = {1, 2, 3}, we know, *NG₁* and *NG₂* are *NG*-groups. We assume *NG₁* ∪ *NG₂* is *NG*- groups. But, it contradiction with example 2-10. Since order of *NG₁* ∪ *NG₂* is four. In addition, *NG₁* ∩ *NG₂*= ϕ, and ϕ ⊂ $S_3$. Therefore, *NG₁* ∩ *NG₂* is not *NG*-group.

## 3- New results

In this section, we introduce the new concepts in NG groups.

**Definition 3-1:** An element *f* of NG is regular if there exists *y*∈ *NG* such that *fyf* =*y*. The NG is regular if all its elements are regular.

Note that NG-groups are regular, since we may choose $y = f^{-1}$.

**Example 3-2**: we consider the example 2-10. If we take NG= $G_1$={(1,1,3),(3,3,1)} ;

If we take *f*=(1,1,3) and *y*=(3,3,1), then *fyf*= (1,1,3) (3,3,1) (1,1,3)= (3,3,1) =*y*.

And if we take *f*= (3,3,1) and *y*=(1,1,3), then *fyf*= (3,3,1) (1,1,3) (3,3,1) = (1,1,3)=*y*.

From definition 3-1, $G_1$ is regular; we can see all NG are regular.

Regularity is equivalent to a condition which appears formally to be stronger:
We can make this a strong condition

**Proposition 3.3** If *f* ∈ *NG* is regular, then there exists *y* ∈ NG such that *fyf* = *y* and *yfy* =*f*.
*Proof*
Choose *f* such that *fyf* = *y*, and set *y* = *yfy*.
Then *fyf* = *fyfyf* = *fyf* = *y*, *yfy* = *yfyfyfy* = *yfyfy* = *yfy* = *f*.

**Definition 3-4:** Suppose that *NG* group, then it is an inverse NG if for each a there is a unique *y* such that *fyf* = *y* and *yfy* = *f*.

The element *y* is called the inverse of *f*.

**Proposition 3.4**: Let NG be an inverse group, each element of NG has a unique inverse.

**Conclusion**: We investigate and present a new result in NG groups such as the definition of regular *NG* and inverse *NG*.

**ACKNOWLEDGMENTS** I would like to thank staff of Higher institute of Suluq for help on this paper.

Rreferences**:**